\documentstyle[12pt]{article}

\newcommand{\beq}{\begin{equation}}
\newcommand{\eeq}{\end{equation}}
\def\square{\hfill{\vrule height6pt width6pt            
depth1pt} \break \vspace{.01cm}}                        

\begin{document}

\rightline{ITEP-TH-46/99}
\bigskip

{\centerline {\bf
Note on the Ruijsenaars-Schneider model.}}

{\centerline {\bf V.Vakulenko}}

  {\footnotesize {\bf Abstract.}
We study vector bundles with some
additional structures
on an elliptic curve and show how there are related to the elliptic
Ruijsenaars-Schneider model.}

  {\bf 0.} Introduction.

  Section {\bf 1} is a brief review of some known facts on
theta-functions, elliptic solutions of the Yang-Baxter equation and
so on.

  In section {\bf 2} we obtain an interpretation of the matrix of
intertwining vectors as a map between some vector bundles on elliptic
curve. It gives the possibility to see that the factorized
Sklyanin-Hasegawa $\bf L$-operator (more precisely, it's "classical"
part) is connected with the certain diagrams. This diagrams are very
similar to that used in works of Drinfeld, where it is called $F$-sheaves.
Although our diagrams doesn't coincide with $F$-sheaf we show that
properties of $F$-sheaves has an analogues in our situation. Also we
discuss how quantum factorized $\bf L$-operator may arise.

  In section {\bf 3} we show that diagrams, studied in this article,
(let us call them $f$-sheaves) are naturally connected with the
Ruijsenaars-Schneider model. We calculate  Backlund transformations
in the RS model and all demanded properties, listed in [KS],
are naturally appear.

  Note that geometric sense of $L$-operators for Calogero-Moser models
was discussed in [ER],[N].

  {\bf 1.}Here we collect some useful facts and formulas. The theta-function
theory are contained in [Fay],[Mu]. For the elliptic $R$-matrix, factorized
$\bf L$-operators see the reference list in [H2].

  The theta-function with characteristics is

\begin{equation}
\theta\left[\begin{array}{cc}a\\b\end{array}\right](z,\tau)=
\sum_{n\in{\bf Z}} \exp{(\pi i(n+a)^{2}\tau+2\pi i(n+a)(z+b))}
\end{equation}
where $a,b$ are rational numbers,$ \tau\in{\bf C},Im\tau>0$.

  For $n\in{\bf Z}_{>0}$ and $j\in{\bf Z}/n{\bf Z}$ we put
\begin{equation}
\theta^{(j)}(z,\tau)= \theta
\left[\begin{array}{cc}\frac{1}{2}-\frac{j}{n}\\0\end{array}\right]
(z+\frac{1}{2},n\tau)
\end{equation}
\begin{equation}
\theta_{j}(z)=
\theta\left[\begin{array}{cc}\frac{1}{2}-\frac{j}{n}\\0\end{array}
\right](n(z+\frac{1}{2}),n\tau)
\end{equation}
and let
\begin{equation}
\theta(z)= \theta\left[\begin{array}{cc}\frac{1}{2}\\
\frac{1}{2}\end{array}\right]
(z,\tau).
\end{equation}
  Then  Belavin's $R$-matrix is [RT]
\begin{equation}
R(z)^{i j}_{i' j'}= \delta_{i+j,i'+j'}\frac{\theta^{(i'-j')}(z+\eta)}
{\theta^{(i'-i)}(\eta)\theta^{(i-j')}(z)}\frac{\prod_{k=0}^{n-1}
\theta ^{(k)}(z)} {\prod_{k=1}^{n-1}\theta^{(k)}(0)}
\end{equation}
where $\eta\in{\bf C}$ is a parameter.

  Let $\epsilon_{i},i=1,...n$ be an orthonormal basis vectors in an
$n$-dimensional vector space: $<\epsilon_{i},\epsilon_{j}>=\delta_{i,j},
{\bf C}^{n}=\oplus_{i=1}^{n}\bf{C}\epsilon_{i}$. Put
$\bar\epsilon_{i}=\epsilon_{i}-\frac{1}{n}\sum_{k=1}^{n}\epsilon_{k}$.
Then for $\lambda\in{\bf C}^{n}$

\begin{equation}
\phi(z)^{\lambda+\eta\bar\epsilon_{k}}_{\lambda ,i}= \theta_{i}(\frac{z}{n}-
<\lambda,\bar\epsilon_{k}>)/\sqrt{-1}{\bf\eta}(\tau)
\end{equation}
are called intertwining vectors [Bax],[JMO]. Here ${\bf\eta}(\tau)$ is the
Dedekind eta-function
$${\bf\eta}(\tau)= \exp(\pi i\tau/12)\prod_{m=1}^{\infty}(1-\exp 2\pi i\tau).$$

  The matrix inverse to $\phi(z)$ is denoted by
$\bar\phi(z)_{\lambda}^{\lambda+\eta\bar\epsilon_{k},j}$.

\begin{equation}
\sum_{i=1}^{n} \phi(z)_{\lambda,i}^{\lambda+\eta\bar\epsilon_{k'}}
\bar\phi(z)_{\lambda}^{\lambda+\eta\bar\epsilon_{k},i}= \delta_{k,k'},
\qquad
\sum_{k=1}^{n}
\phi(z)_{\lambda,i}^{\lambda+\eta\bar\epsilon_{k}}
\bar\phi(z)_{\lambda}^{\lambda+\eta\bar\epsilon_{k},i'}
= \delta_{i,i'}.
\end{equation}

  There is the useful formula
$$
\sum_{i=1}^{n}
\bar\phi(z)_{\mu}^{\mu+\eta\bar\epsilon_{k},i}
\phi(z+u)_{\lambda,i}^{\lambda+\eta\bar\epsilon_{k'}}=
$$
\beq
\label{f-f}
=\frac{\theta(z+\frac{u}{n}+<\mu,\bar\epsilon_{k}>-
<\lambda,\bar\epsilon_{k'}>)}
{\theta(z)}
\prod_{l\neq k}
\frac{\theta(\frac{u}{n}+<\mu,\bar\epsilon_{l}>-
<\lambda,\bar\epsilon_{k'}>)}
{\theta(<\mu,\bar\epsilon_{l}>-
<\mu,\bar\epsilon_{k}>)}
\end{equation}
which follows from the determinant formula
\begin{equation}
\det(\theta_{i}(z_{j}))_{i,j=1,...n}=
(-1)^{n-1}\frac{\theta(\sum_{j}z_{j})}{\sqrt{-1}{\bf\eta}(\tau)}
\prod_{i<j}
\frac{\theta(z_{j}-z_{i})}
{\sqrt{-1}{\bf\eta}(\tau)}
\end{equation}
  Note also the following property of the matrix $\phi(z)$:

\begin{equation}
\label{detphi}
\det(\phi(z)_{\lambda,i}^{\lambda+\bar\epsilon_{k}})_{i,k=1,...n}=
(-1)^{n-1}\frac{\theta(z)}{\sqrt{-1}{\bf\eta}(\tau)}
\prod_{i<j}
\frac{\theta(\lambda_{i}-\lambda_{j})}
{\sqrt{-1}{\bf\eta}(\tau)}
\end{equation}

  In [S] for $n=2$ and in [H1] for any $n$ ${\bf L}$-operator intertwined by
Belavin's $R$-matrix was found
\begin{equation}
\label{L}
{\bf L}(z)^{j}_{i}= \sum_{k=1}^{n}
\bar\phi(z)_{\lambda}^{\lambda+\eta\bar\epsilon_{k},j}
\phi(z+\sigma\eta)_{\lambda,i}^{\lambda+\eta\bar\epsilon_{k}}\bar T_{k}
\end{equation}
where $\lambda\in{\bf C}^{n},\sigma\in{\bf C}$ is arbitrary and
$$
({\bar T}_{k}f)(\lambda)= f(\lambda+\eta\bar\epsilon_{k})
$$

  Also we need some properties of the function $\theta(z)$. Let
\begin{equation}
\Phi_{z}(x)=
\frac{\theta(z+x)}{\theta(z)\theta(x)}
\qquad
\zeta(z)=\frac{\theta'(z)}{\theta(z)}
\end{equation}
  The functional relation for $\Phi_{z}(x)$ is
\begin{equation}
\label{a}
\Phi_{z}(x)
\Phi_{z}(y)=
\Phi_{z}(x+y)(\zeta(z)+\zeta(x)+\zeta(y)-\zeta(z+x+y))
\end{equation}
  An elliptic form of the Lagrange interpolation formula is
\begin{equation}
\label{b}
\prod_{i=1}^{N}
\frac{\theta(z-x_{i})}
{\theta(z-y_{i})}=
\sum_{i=1}^{N}
(\zeta(z-y_{i})-\zeta(x-y_{i}))
\frac{\prod_{j=1}^{N}
\theta(y_{i}-x_{j})}
{\prod_{j\neq i}
\theta(y_{ij})}
\end{equation}
when
$\sum_{i=1}^{N}(x_{i}-y_{i})= 0$.
Here $x$ is any of $x_{i}$.

  Note also the useful relation
\begin{equation}
\label{c}
\sum_{i=1}^{N}
\frac{\prod_{j=1}^{N}
\theta(y_{i}-x_{j})}
{\prod_{j\neq i}
\theta(y_{ij})}
=0
\qquad
\mbox{when}
\qquad
\sum_{i=1}^{N}(x_{i}-y_{i})= 0.
\end{equation}

  {\bf 2.}Let $X$ be an elliptic curve with periods $1,\tau$.
Here we recall some general facts about vector bundles on elliptic curve.

  The description of the set of semistable bundles of rank $r$ and degree
$d$ on an elliptic curve depends crucially on arithmetic properties of
$r$ and $d$. If $r$ and $d$ are coprime $(r,d)=1$, then such semistable
bundles are parametrized by the elliptic curve itself. Then $(r,d)=0$
the moduli space of semistable bundles is isomorphic to $S^{r}X$-$r$-th
symmetric power of $X$. This cases will be most important for us.
For example, generic vector bundle of rank $n$ and degree zero is
equivalent to $\oplus_{i=1}^{n}{\cal O}_{X}(P_{i}-P_{0})$, where
$P_{0},P_{i}\in X, i=1,...n$ and we fix $P_{0}$ by $z(P_{0})= 0$.

  An elementary modification of the bundle $\cal F$ at the point $P\in X$
is the vector bundle $\cal E$ of the same rank (but another degree) which
is isomorphic to $\cal F$ everywhere except the point $P$. Namely there
is the exact sequence
$$
0\longrightarrow
{\cal F}
\longrightarrow
{\cal E}
\longrightarrow
{\cal Q}
\longrightarrow 0
$$
where ${\cal Q}$ is the skyscraper sheaf with support at the point $P$.
If the kernel of the map
${\cal F}|_{P}\longrightarrow {\cal E}|_{P}$
is one-dimensional the determinant of the map
${\cal F}\longrightarrow {\cal E}$
has precisely one zero on the curve $X$. Due to (\ref{detphi}) the
determinant of $\phi(z)$ also has precisely one zero and this is not
accidential. We show that the matrix of intertwining vectors coincide
with the matrix of an elementary modification of a some vector bundle.
We explicitely describe this modification using Fourier transform [Muk].

Namely, let $p_{1},p_{2}$ are corresponding projections
$X\stackrel{p_{1}}{\leftarrow}
X\times X
\stackrel{p_{2}}{\rightarrow} X$ to
the first and second factors and $\cal{P}$ is the Poincare line bundle on
$X\times X$
$${\cal P}=
p_{1}^{*}{\cal O}_{X}(-P_{0})\otimes p_{2}^{*}{\cal O}_{X}(-P_{0})\otimes
{\cal O}_{X\times X}
(\Delta),
$$
$\Delta\subset X\times X$ is the diagonal. Then if ${\cal F}$ is a
coherent sheaf on $X$ it's image under Fourier transform is by definition
$R{\cal T} ({\cal F})= Rp_{2*}p_{1}^{*}({\cal P}\otimes{\cal F})$.

  Let $P_{\lambda_{1}},...P_{\lambda_{n}}$ be the points of $X$,
$z(P_{\lambda_{i}})=\lambda_{i}, i=1,...n$ and ${\cal L}_{1,n}= {\cal O}_{X}
(P_{\lambda_{1}}+...+P_{\lambda_{n}})$ is the line bundle of degree $n$.
Consider the exact sequence
$$
0\longrightarrow{\cal O}_{X}\longrightarrow {\cal L}_{1,n}
\stackrel{ev}{\longrightarrow}
\oplus_{i=1}^{n}{\cal L}_{1,n}|_{P_{\lambda_{i}}}
\longrightarrow 0
$$
and applying to them Fourier
transform we get

$$ 0\longrightarrow p_{2*}p_{1}^{*}({\cal L}_{1,n}\otimes{\cal P})
\stackrel{i_{\alpha}}{\longrightarrow}
\oplus_{i=1}^{n}p_{2*}p_{1}^{*}
({\cal L}_{1,n}|_{P_{\lambda_{i}}})
\longrightarrow R^{1}p_{2*}{\cal P}
\longrightarrow 0
$$

  The sheaf $ R^{1}p_{2*}{\cal P}$ is a scyscraper sheaf with support
at the point $P_{0}$, ${\cal F}_{n,-1}
= p_{2*}p_{1}^{*}({\cal L}_{1,n}\otimes{\cal P})$
and ${\cal E}_{n,0}=
\oplus_{i=1}^{n}p_{2*}p_{1}^{*} ({\cal
L}_{1,n}|_{P_{\lambda_{i}}})$ 
are vector bundles of rank $n$ and degrees
$-1$
and $0$. Note that ${\cal E}_{n,0}=
\oplus_{i=1}^{n}{\cal O}
(P_{\lambda_{i}}-P_{0})$
and the fiber of the bundle ${\cal F}_{n,-1}$ over the point $P_{z}$ is
\begin{equation}
{\cal F}_{n,-1}|_{P_{z}}= \oplus_{i=1}^{n}{\bf C}\theta_{i}
(\frac{\sum_{j=1}^{n}\lambda_{j}+z}{n}-x)
\end{equation}
where $x\in X_{1}, z\in X_{2}$.

  The map $i_{\alpha}$ is simply the the sum of the
evaluations at the points
$P_{\lambda_{k}}$ and sends a section $s(x)= \sum_{i=1}^{n}a_{i}
\theta_{i}(\frac{\sum \lambda_{j}+z}{n}-x)$ to the column
$(s(\lambda_{1}),...s(\lambda_{n}))^{t}$. The kernel of $i_{\alpha}$
over the point $P_{0}$ is generated by $s_{\lambda}(x)= \prod_{j=1}^{n}
\theta(\lambda_{j}-x)$.

  Then in a given basises in fibers of bundles ${\cal F}_{n,-1},
{\cal E}_{n,0}$
the matrix elements of $i_{\alpha}$ are
\begin{equation}
\label{ev}
{(i_{\alpha})}_{i}^{k}(z)=
\theta_{i}(\frac{\sum\lambda_{j}+z}{n}-\lambda_{k})=
\theta_{i}(\frac{z}{n}-<\lambda,{\bar\epsilon}_{k}>)
\end{equation}
which up to multiple coincides with $\phi(z)_{\lambda,i}^{\lambda+
\eta\bar\epsilon_{k}}$.

  Next, let us consider the Sklyanin-Hasegawa ${\bf L}$-operator
(\ref{L})(and put
$\sigma=1$):
\begin{equation}
\label{L1}
{\bf L}(z)^{j}_{i}= \sum_{k=1}^{n}
\bar\phi(z)_{\lambda}^{\lambda+\eta\bar\epsilon_{k},j}
\phi(z+\eta)_{\lambda,i}^{\lambda+\eta\bar\epsilon_{k}}\bar T_{k}
\end{equation}

  From the previous discussion it follows that the "classical part"
of the $\bf L$-operator (in indices $i,j$)
is described by  a diagram
\beq
\label{rdiagr}
{\cal F}
\stackrel{i_{\alpha}}
{\longrightarrow}
{\cal E}
\stackrel{i_{\beta}}
{\longleftarrow}
{\cal F}\otimes {\bf \xi}
\eeq
where $i_{\alpha}$ is an elementary modifications at the point $P_{-\eta}$,
$i_{\beta}$ is an elementary modification at the point $P_{0}$,
$\bf{\cal \xi}$ is a line bundle of degree zero, ${\bf{\cal \xi}}=
{\cal O}_{X}(P_{0}-P_{\eta/n})$. Namely, "classical part"
of $L$-operator is a rational map from ${\cal F}$ to ${\cal F}\otimes
{\bf{\cal \xi}}$ with zero at $P_{-\eta}$ and pole at $P_{0}$.

  A more general diagrams of this kind appears in the works of Drinfeld [Dr]
in the study of the Langlands correspondence over functional field. Those
diagrams are called $F$-sheaves. This is the motivation for the
following

  {\bf Definition.} Right $f$-sheaf of rank $n$ is the diagram
\beq
\label{rf}
\begin{array}{ccc}{\cal F}&{}&{}
\\{}&\stackrel{i_{\alpha}}{\searrow} &{}
\\{}&{}&{\cal E}\\
{}&\stackrel{i_{\beta}}{\nearrow}&{}\\
{\cal F}\otimes{\bf \xi}&{}&{}
\end{array}
\eeq
where ${\cal F},{\cal E}$ are vector bundles of rank $n$, ${\bf \xi}$ is
the line bundle of degree zero, $i_{\alpha},i_{\beta}$ are injective
and cokernels of $i_{\alpha},i_{\beta}$ are one-dimensional and lies over
the points $P_{\alpha},P_{\beta}\in X, P_{\alpha}\neq P_{\beta}$.

  Left $f$-sheaf of rank $n$ is the diagram
\beq
\label{lf}
\begin{array}{ccc}{}&{}&{\cal F}
\\{}&\stackrel{i_{\beta}'}{\nearrow}&{}\\
{\cal G}&{}&{}\\
{}&\stackrel{i_{\alpha}'}{\searrow} &{}
\\{}&{}&{\cal F}\otimes{\bf \xi}
\end{array}
\eeq
where ${\cal F},{\cal G}$ are vector bundles of rank $n$, ${\bf \xi}$ is
the line bundle of degree zero, $i_{\alpha}',i_{\beta}'$ are injective
and cokernel s of $i_{\alpha}',i_{\beta}'$ are one-dimensional over
the points $P_{\alpha},P_{\beta}\in X, P_{\alpha}\neq P_{\beta}$.
The (classical) $L$-operator, associated to the $f$-sheaf is a
corresponding rational map from $\cal F$ to ${\cal F}\otimes \xi$.

  Actually the right and left $f$-sheaves are the same. Changing the
order of elementary modifications we move from (\ref{rf}) to (\ref{lf})
and vice versa.

  {\bf Example.} Let ${\cal E}$ be the vector bundle of rank $n$ and degree
zero, $P_{\alpha}= P_{v},P_{\beta}= P_{v+\eta}$ and
\beq
\begin{array}{ccccc}{}&{}&{\cal E}&{}&{}\\
{}&\stackrel{i_{\beta}}
{\nearrow}&{}&
\stackrel{i_{\alpha}}
{\searrow}&{}\\
{\cal F}&{}&{}&{}&
{\cal G}\\
{}&\stackrel{i_{\alpha}'}
{\searrow}&{}
&\stackrel{i_{\beta}}
{\nearrow}&{}\\
{}&{}&
{\cal E}\otimes {\bf \xi}&
{}&{}
\end{array}
\eeq
i.e. ${\cal E}$ is the right and the left $f$-sheaf simultaneousely.
Then the commutativity of the elementary modifications is illustrated
by the identity:
\beq
\label{commute}
\sum_{i}
\phi(z-v)_{\lambda,i}^{\lambda+\eta{\bar\epsilon}_{k'}}
{\bar\phi}(z-v-\eta)_{\lambda}^{\lambda+\eta{\bar\epsilon}_{k},i}=
\frac{\prod_{m\neq k'}\theta(\lambda_{k'm})}
{\prod_{m\neq k}\theta(\lambda_{mk})}
\prod_{l}
\frac{\theta(\lambda_{lk'}+\frac{\eta}{n})}
{\theta(\lambda_{kl}+\frac{\eta}{n})}\times
\eeq
$$
\times
\sum_{i}
{\bar\phi}(z-v-\eta)_{-\lambda}^{-\lambda+\eta{\bar\epsilon}_{k'},i}
\phi(z-v)_{-\lambda,i}^{-\lambda+\eta{\bar\epsilon}_{k}}.
$$
  The appearence of $-\lambda$ in the intertwining vectors follows
from the fact that fiberwise the map ${\cal E}
\stackrel{i_{\alpha}}{\longrightarrow}
{\cal G}$ is dual to the evaluations map (\ref{ev}). Note also that
in trigonometric case the equation (\ref{commute}) appears in [BKMS]
(equation (3.5a)).

  There are exist another way to transform the right $f$-sheaf
to the left and vice versa. Namely if (\ref{rf}) is given,
the diagram
$
{\cal E}
\stackrel{i_{\beta}}
{\longleftarrow}
{\cal F}\otimes{\bf \xi}
\stackrel{i_{\alpha}}
{\longrightarrow}
{\cal E}\otimes {\bf \xi}
$
is the left $f$-sheaf and similary from (\ref{lf}) we can obtain the
diagram
$
{\cal G}
\stackrel{i_{\alpha}'}
{\longrightarrow}
{\cal F}\otimes {\bf \xi}
\stackrel{i_{\beta}'}
{\longleftarrow}
{\cal G}\otimes {\bf \xi}.
$

  {\bf Example.} Let $\cal F$ be a right $f$-sheaf of rank $n$. Then the
transformation of (\ref{rdiagr}) to the diagram
${\cal E}\otimes\xi ^{-1} \
\longleftarrow
{\cal F}
\longrightarrow
{\cal E}
$
change the $\bf L$-operator (\ref{L1}) to
\beq
\label{conjl}
{\bf L}(z)_{i}^{j}\longrightarrow \sum_{i,j}
\phi(z)^{\lambda+\eta\bar\epsilon_{k^{,}}}_{\lambda ,j}
\bar\phi(z)_{\lambda}^{\lambda+\eta\bar\epsilon_{k},i}
{\bf L}(z)_{i}^{j}=
\frac{\theta(z+\frac{\eta}{n}+\lambda_{kk'})}{\theta(z)}
\prod_{j\neq k}\frac{\theta(\lambda_{jk'}+\frac{\eta}{n})}
{\theta(\lambda_{jk})}
\bar T_{k^{,}}
\eeq
This conjugation of the $L$-operator (\ref{L1}) was noted in [H2].

  Using this operations we can extend (\ref{rf}) or (\ref{lf}) to the
infinite diagram:
\beq
\label{inf}
\begin{array}{ccccccc}{}&{}&{}&{\cal F}_{0}\otimes \xi^{-1}&
{}&{}&{}\\
\searrow &{}&\nearrow &{}&\searrow &{}&\nearrow\\
{}&{\cal F}_{-1}&{}&{}&{}&{\cal F}_{1}\otimes\xi^{-1}&{}\\
\nearrow &{}&\searrow &{}&\nearrow &{}&\searrow \\
{}&{}&{}&{\cal F}_{0}&{}&{}&{}\\
\searrow &{}&\nearrow &{}&\searrow &{}&\nearrow\\
{}&{\cal F}_{-1}\otimes \xi&{}&{}&{}&{\cal F}_{1}&{}\\
\nearrow &{}&\searrow &{}&\nearrow &{}&\searrow\\
{}&{}&{}&{\cal F}_{0}\otimes \xi&{}&{}&{}
\end{array}
\eeq
where the maps in the NW-SE direction are elementary modifications at the
point $P_{\alpha}$ and the maps in SW-NE direction are elementary
modifications at the point $P_{\beta}$.
Actually the diagram (\ref{inf}) is defined by any $f$-sheaf contained
in it.
It is convinient to parametrize the diagrams (\ref{inf}) by it's
$f$-sheaves of degree $-1$ or zero. This is because the description
of generic vector bundles of degree $-1$ or zero is very simple.

  The diagrams (\ref{inf}) are useful when we consider the action of
(some) elementary modifications on $f$-sheaves. Let ${\cal F}
\longrightarrow {\cal F}'$ be the simplest elementary modification
of the bundle $\cal F$ at the point $P\neq P_{\alpha},P_{\beta}$
with $deg{\cal F}'= deg{\cal F}+1$. The diagram (\ref{rf}) gives the
identification of fibers
${\cal F}|_{P}
\to {\cal E}|_{P}\to ({\cal F}\otimes \xi)|_{P}
$
and therefore we can apply the elementary modification at the point $P$
to the bundles
${\cal E}$ and ${\cal F}\otimes \xi$. But in general the $f$-sheaf
(\ref{rf}) doesn't transform to another $f$-sheaf. It is possible iff
the kernel of the map
${\cal F}|_{P} \to
{\cal F}'|_{P} $
is the eigenvector of the composition
\beq
\label{composit}
{\cal F}|_{P}\to{\cal E}|_{P}\to({\cal F}\otimes \xi)|_{P}
\to{\cal F}|_{P}
\eeq
where first and second arrows comes from the $f$-sheaf structure and
the last arrow is a natural isomorphism. In other words the kernel of
the elementary modification must be an eigenvector of the classical
$L$-operator at the point $P$.

  Althought under the elementary modification $f$-sheaf transforms
to $f$-sheafs of another degree the whole diagram (\ref{inf}) is mapped to
a diagram of the same kind and this is important for explicit formulas of
the action of the elementary modifications. In the next paragraph we
show that this action coincides with Backlund transformations of the
Ruijsenaars-Schneider (RS) model [KS].




  Now we consider (\ref{rdiagr}) more closely and show that it is connected
not only with $L$-operator of the classical RS model but also with
"quantum" Sklyanin-Hasegawa $\bf L$-operator.

  To define the maps in the diagram we need the one-dimensional subspace in
$$
{\cal F}|_{P_{-\eta}}= \sum_{i=1}^{n}{\bf C}\theta_{i}(\frac{\Lambda}{n}-x)
$$
and $(n-1)$-dimensional subspace in
$$
{\cal F}|_{P_{0}}= \sum_{i=1}^{n}
{\bf C}\theta_{i}(\frac{\Lambda+\eta}{n}-x)
$$
(or equivalently a
linear functional on ${\cal F}|_{P_{0}}$) where we set $\Lambda=
\sum_{j}\lambda_{j}$.
Let $s_{\lambda}(x)=\prod
_{l=1}^{n}\theta(\lambda_{l}-x)\in{\cal F}|_{P_{-\eta}}$ lies in the kernel
of $i_{\alpha}$. Then the classical $L$-operator takes the form:
\beq
\label{str}
L(z)_{i}^{j}=
\sum_{k=1}^{n}
{\bar\phi}(z)_{\lambda}^{\lambda+\eta{\bar\epsilon}_{k},j}
\phi(z+\eta)_{\lambda,i}^{\lambda+\eta{\bar\epsilon}_{k}}
t_{k}
\eeq
where $t_{k},k=1,...n$ are parameters. We need to know how $t_{k}$ are
connected with the cokernel of $i_{\beta}$. To do this we compute
the linear functional $\psi_{t}$ on ${\cal F}|_{P_{0}}$ which
arise from
the $L(z)$ in the form (\ref{str})
$$
\theta_{i}(\frac{\Lambda+\eta}{n}-x)/{\sqrt -1}{\bf\eta}(\tau)\longrightarrow
\sum_{k,j}
\phi(\eta)_{\lambda,i}^{\lambda+\eta\bar\epsilon_{k}}t_{k}
{\tilde \phi}(0)_{\lambda}^{\lambda+\eta\bar\epsilon_{k},j}
\theta_{j}
(\frac{\Lambda}{n}-x)=
$$
$$
=\sum_{k}
\phi(\eta)_{\lambda,i}^{\lambda+\eta\bar\epsilon_{k}}
t_{k}\theta(<\lambda,{\bar\epsilon}_{k}>+\frac{\Lambda}{n}-x)\prod_{l\neq k}
\frac{\theta(<\lambda,{\bar\epsilon}_{l}>+\frac{\Lambda}{n}-x)}{\theta
(\lambda_{lk})}=
$$
$$
=\sum_{k}\frac{\phi(\eta)_{\lambda,i}^{\lambda+\eta\bar\epsilon_{k}}t_{k}}
{\prod_{l\neq k}\theta(\lambda_{lk})}\times\prod_{j=1}^{n}\theta
(\lambda_{j}-x)
$$
where
\beq
\label{tildph}
{\tilde \phi}(0)_{\lambda}^{\lambda+\eta\bar\epsilon_{k},j}=
(\theta(z)
{\bar\phi}(z)_{\lambda}^{\lambda+\eta{\bar\epsilon}_{k},j})|_{z=0}.
\eeq
Therefore
$$
<\psi_{t},\theta_{i}(\frac{\Lambda+\eta}{n}-x)/{\sqrt -1}{\bf \eta}(\tau)>=
\sum_{k}\frac{t_{k}}
{\prod_{l\neq k}\theta(\lambda_{lk})}
\phi(\eta)_{\lambda,i}^{\lambda+\eta\bar\epsilon_{k}}=
\sum_{k}\frac{t_{k}}{\prod_{l\neq k}\theta(\lambda_{lk})}
\theta_{i}(\frac{\Lambda+\eta}{n}-\lambda_{k})
$$
or
$$
t_{k}=\prod_{l\neq k}\theta(\lambda_{lk})\sum_{m}
{\bar\phi}(\eta)_{\lambda}^{\lambda+\eta\bar\epsilon_{k},m}
<\psi_{t},\theta_{m}(\frac{\Lambda+\eta}{n}-x)>
$$
  Over the point $P_{z}$ the $L$-operator (\ref{str}) acts from
$$
{\cal F}|_{P_{z}}
= \sum_{i=1}^{n}{\bf C}\theta_{i}(\frac{\Lambda+z+\eta}{n}-x)
$$
to
$$({\cal F}\otimes\xi)|_{P_{z}}=
\sum_{i=1}^{n}{\bf C}\theta_{i}(\frac{\Lambda+z}{n}-x).
$$
"Define"
${\hat L}(z)$ as a map from
${\cal F}|_{P_{z}}\otimes
{\cal F}|_{P_{-\eta}}$ to
$({\cal F}\otimes\xi)|_{P_{z}}\otimes
{\cal F}|_{P_{0}}$ such that
$$
<\psi_{t},{\hat L}(z)
(\theta_{i}(\frac{\Lambda+z+\eta}{n}-x)
\otimes
s_{\lambda})>=
\sum_{j}L_{i}^{j}(z)
\theta_{j}(\frac{\Lambda+z}{n}-x)=
$$
$$
=\sum_{k,j}
{\bar\phi}(z)_{\lambda}
^{\lambda+\eta{\bar\epsilon}_{k},j}
\phi(z+\eta)_{\lambda,i}^{\lambda+\eta\bar\epsilon_{k}}
t_{k}
\theta_{j}(\frac{\Lambda+z}{n}-x).
$$
Then
$$
<\psi_{t},{\hat L}(z)
(\theta_{i}(\frac{\Lambda+z+\eta}{n}-x)
\otimes
s_{\lambda})>=
\sum_{k,j}
{\bar\phi}(z)_{\lambda}^{\lambda+\eta{\bar\epsilon}_{k},j}
\phi(z+\eta)_{\lambda,i}^{\lambda+\eta\bar\epsilon_{k}}
\prod_{l\neq k}\theta(\lambda_{lk})\times
$$
$$
\times
\sum_{m}
{\bar\phi}(\eta)_{\lambda}^{\lambda+\eta{\bar\epsilon}_{k},m}
<\psi_{t},\theta_{m}(\frac{\Lambda+\eta}{n}-x)>
\theta_{j}(\frac{\Lambda+z}{n}-x)=
$$
$$
=\sum_{k}
{\bar\phi}(z)_{\lambda}^{\lambda+\eta{\bar\epsilon}_{k},j}
\phi(z+\eta)_{\lambda,i}^{\lambda+\eta\bar\epsilon_{k}}
\theta_{j}(\frac{\Lambda+z}{n}-x)
<\psi_{t},\theta(\eta+\lambda_{k}-x)\prod_{l\neq k}
\theta(\lambda_{l}-x)>
\frac
{{\sqrt -1}{\bf \eta}(\tau)}
{\theta(\eta)}
$$
where we use the identity
$$
\sum_{m}
{\bar\phi}(\eta)_{\lambda}^{\lambda+\eta{\bar\epsilon}_{k},m}
\theta_{m}(\frac{\Lambda+\eta}{n}-x)={\sqrt -1}{\bf \eta}(\tau)
\frac{\theta(\eta+\lambda_{k}-x)}{\theta(\eta)}\prod_{l\neq k}
\frac{\theta(\lambda_{l}-x)}{\theta(\lambda_{lk})}.
$$
Therefore
$$
{\hat L}(z)_{i}^{j}=\frac{{\sqrt -1}{\bf \eta}(\tau)}{\theta(\eta)}
\sum_{k}
{\bar\phi}(z)_{\lambda}^{\lambda+\eta{\bar\epsilon}_{k},j}
\phi(z+\eta)_{\lambda,i}^{\lambda+\eta\bar\epsilon_{k}}
T_{k}
$$
where $(T_{k}f)(\lambda_{1},...\lambda_{n})= f(\lambda_{1},...,
\lambda_{k}+\eta,
...\lambda_{n})$. Then $f(\lambda)$ is independent on the sum
$\sum \lambda_{i}$,$ T_{k}f= {\bar T}_{k}f$ and the formula up to
multiple coincides with (\ref{L1}). But this "definition" is incorrect
and this point deserves the further study.

  {\bf 3.}In this paragraph we study the connection of diagrams (\ref{rf}),
(\ref{lf})
with RS model. In [H2] it was pointed that ${\bf L}$-operator (\ref{conjl})
gives Lax matrix for RS system and here we obtain  Backlund
transformations [KS] and discrete time dinamics [NRK] of this Lax
matrix.

  Let
\beq
\label{L3}
L(z)_{i}^{j}=
\sum_{k=1}^{n}
{\bar\phi}(z-v-\eta)_{\lambda}^{\lambda+\eta{\bar\epsilon}_{k},j}
\phi(z-v)_{\lambda,i}^{\lambda+\eta{\bar\epsilon}_{k}}
t_{k}
\eeq
be the $L$-operator corresponding to the diagram (\ref{rf}) with zero at the
point $P_{v}$ and pole at the point $P_{v+\eta}$ and $deg {\cal F}=-1$.
An elementary
modification ${\cal F}\stackrel{\phi_{u}}{\longrightarrow}{\cal G}$
of the bundle ${\cal F}={\cal F}_{-1}$
at the point $P_{u}\neq P_{v},P_{v+\eta}$
transforms the diagram (\ref{rf}) into a diagram of the
same kind iff the kernel of the map $\phi_{u}$ over the point $P_{u}$ is an
eigenvector of the $L$-operator (\ref{L3}). In this case "new" $L$-operator
has the same properties as the "old" one: it is invertible everywhere
except the points $P_{v},P_{v+\eta}$ with zero at $P_{v}$ and pole
at $P_{v+\eta}$. Consider this more closely.

  Let ${\cal G}= \oplus_{i=1}^{n}{\cal O}_{X}(P_{\mu_{i}}-P_{0})$
and the kernel of $\phi_{u}$ is
$ {\bf C}s_{\mu}(x)=
{\bf C}\prod_{s}\theta(x-\mu_{s})
\in{\cal F}|_{P_{u}}$. Then this vector is an eigenvector of $L(u)$
with eigenvalue $e^{c}$
iff
$t_{k}s_{\mu}(\lambda_{k})= e^{c}s_{\mu}(\lambda_{k}+\frac{\eta}{n})$
for any $k=1,...n$, as follows from the explicit description of the
maps in the diagram (\ref{rf}).
Therefore
\beq
\label{t}
t_{k}= e^{c}\prod_{s}
\frac{\theta(\lambda_{k}-\mu_{s}+\frac{\eta}{n})}
{\theta(\lambda_{k}-\mu_{s})}
\eeq

  As it was noted previousely we can think about the elementary modification
$\phi_{u}$ as the map from the diagram (\ref{inf}) to another such
diagram, constructed from a new $f$-sheaf. Let us consider a part of this
map:
\beq
\label{mapff'}
\begin{array}{ccccccc}{}&{}&{\cal F}'_{0}\otimes{\bf \xi}^{-1}={\cal G}&
\stackrel{\phi_{u}}{\longleftarrow}&{\cal F}_{-1}&{}&{}\\
{}&\nearrow &{}&{}&{}&
\stackrel{i_{\alpha}}{\searrow}&{}\\
{\cal F}'_{-1}&{}&{}&{}&{}&{}&{\cal F}_{0}\\
{}&\searrow &{}&{}&{}&
\stackrel{i_{\beta}}{\nearrow}&{}\\
{}&{}&{\cal F}'_{0}&
\stackrel{\phi_{u}}{\longleftarrow}{\cal F}_{-1}
\otimes{\bf \xi}&{}&{}\\
{}&\nearrow &{}&{}&{}&{}&{}\\
{\cal F}'_{-1}\otimes {\bf \xi}&{}&{}&{}&{}&{}&{}
\end{array}
\eeq
  By definition, $L(z)$ is the map from ${\cal F}_{-1}$ to ${\cal F}_{-1}
\otimes {\bf\xi}$, new $L$-operator ${\tilde L}(z)$ is the map from
${\cal F}'_{-1}$ to ${\cal F}'_{-1}\otimes{\bf\xi}$. Let $M(z)$ denote
the map from ${\cal F}_{-1}$ to ${\cal F}'_{-1}$. Then we have the
commutative diagram:
\beq
\label{LLMM}
\begin{array}{ccccc}{}&{\cal F}_{-1}&
\stackrel{L(z)}{\longrightarrow}&{\cal F}_{-1}\otimes {\bf\xi}&{}\\
M(z)&\downarrow &{}&\downarrow &M(z)\\
{}&{\cal F}'_{-1}&
\stackrel{{\tilde L}(z)}{\longrightarrow}&{\cal F}'_{-1}\otimes
{\bf\xi}&{}
\end{array}
\eeq
The vertical arrows are the same due to $ Hom({\cal F}_{-1},
{\cal F}'_{-1})= Hom({\cal F}\otimes {\bf\xi},{\cal F}'_{-1}
\otimes{\bf \xi})$.

  The diagram (\ref{LLMM}) is the discrete Lax equation:
\beq
\label{Laxeq}
\sum_{j}
M(z)_{j}^{j'}L(z)_{i'}^{j}=
\sum_{i}
{\tilde L}(z)_{i}^{j'}M(z)_{i'}^{i}
\eeq
and from the diagram (\ref{mapff'}) we know that
\beq
\label{tildL}
{\tilde L}(z)_{i}^{j}=
\sum_{k}
{\bar\phi}(z-v-\eta)_{\mu}^{\mu+\eta{\bar\epsilon}_{k},j}
\phi(z-v)_{\mu,i}^{\mu+\eta{\bar\epsilon}_{k}}
{\tilde t}_{k}
\eeq
\beq
M(z)_{i}^{j}= \sum_{k}
{\bar\phi}(z-v-\eta)_{\mu}^{\mu+\eta{\bar\epsilon}_{k},j}
\phi(z-u)_{\mu,i}^{\mu+\eta{\bar\epsilon}_{k}}
C_{k}
\eeq
with some ${\tilde t}_{k},C_{k},k=1,...n$. We can use the Lax equation
to express ${\tilde t}_{k},C_{k}$ as functions of $\lambda,\mu$. In this
way it was done in [NRK]. But it is instructive to see how to compute, for
example, $C_{k}$ directly.

  From (\ref{mapff'}) we have the sequence of isomorphisms
$$
{\cal F}_{-1}|_{P_{v}}
\stackrel{\phi_{u}}{\longrightarrow}
{\cal G}|_{P_{v}}
\longrightarrow{\cal F}'_{-1}|_{P_{v}}
$$
It's composition must send the kernel of
${\cal F}_{-1}|_{P_{v}}
\longrightarrow
{\cal F}_{0}|_{P_{v}}$, which is
${\bf C}\prod_{s}\theta(\lambda_{s}-x)$
, to the kernel of
${\cal F}'_{-1}|_{P_{v}}
\longrightarrow
{\cal F}'_{0}|_{P_{v}}$, which is
${\bf C}\prod_{s}\theta(-\frac{\eta}{n}+\mu_{s}-x)$.
Hence, we can normalize $C_{k}$ such that
\beq
\sum_{k,i}
(\prod_{s}\theta(\lambda_{s}-\mu_{k}))
C_{k}
{\bar\phi}(-\eta)_{\mu}^{\mu+\eta{\bar\epsilon}_{k},i}
\theta_{i}(\frac{-\eta+\sum_{j}\mu_{j}}{n}-x)=
{\sqrt -1}{\bf\eta}(\tau)\prod_{s}
\theta(-\frac{\eta}{n}+\mu_{s}-x).
\eeq
The right hand side is equal
(${\tilde \phi}(0)$
was defined in
(\ref{tildph})):
\beq
{\sqrt -1}{\bf\eta}(\tau)\prod_{s}
\theta(-\frac{\eta}{n}+\mu_{s}-x)=
\prod_{l\neq k'}\theta(\lambda_{lk'})
\sum_{i}
{\tilde \phi}(0)_{\mu}^{\mu+\eta{\bar\epsilon}_{k'},i}
\theta_{i}(\frac{-\eta+\sum_{j}\mu_{j}}{n}-x)
\eeq
for any $k'$.
Therefore
\beq
\sum_{k}
(\prod_{s}\theta(\lambda_{s}-\mu_{k}))
C_{k}
{\bar\phi}(-\eta)_{\mu}^{\mu+\eta{\bar\epsilon}_{k},i}=
\prod_{l\neq k'}\theta(\lambda_{lk'})
{\tilde \phi}(0)_{\mu}^{\mu+\eta{\bar\epsilon}_{k'},i}.
\eeq
Multiplying both sides by
$\phi(-\eta)_{\mu,i}^{\mu+\eta{\bar\epsilon}_{k}}$
and summing over $i$ we get
\beq
(\prod_{s}\theta(\lambda_{s}-\mu_{k}))
C_{k}=
\prod_{l\neq k'}\theta(\lambda_{lk'})
\sum_{i}
{\tilde \phi}(0)_{\mu}^{\mu+\eta{\bar\epsilon}_{k'},i}
\phi(-\eta)_{\mu,i}^{\mu+\eta{\bar\epsilon}_{k}}
=\prod_{s}\theta(\mu_{sk}-\frac{\eta}{n})
\eeq
and
\beq
\label{CK}
C_{k}= \prod_{s}\frac{\theta(\mu_{sk}-\frac{\eta}{n})}
{\theta(\lambda_{s}-\mu_{k})}.
\eeq

  Let us check that equations (\ref{t}) and (\ref{CK}) doesn't contradict
to the Lax equation. Consider (\ref{Laxeq}) in detail:
$$
\sum_{k,k'}
{\bar\phi}(z-v-\eta)_{\mu}^{\mu+\eta{\bar\epsilon}_{k'},j'}
(\sum_{j}
\phi(z-u)_{\mu,j}^{\mu+\eta{\bar\epsilon}_{k'}}
{\bar\phi}(z-v-\eta)_{\lambda}^{\lambda+\eta{\bar\epsilon}_{k},j})
\phi(z-v)_{\lambda,i'}^{\lambda+\eta{\bar\epsilon}_{k}}
C_{k'}t_{k}=
$$
\beq
=\sum_{k,k'}
{\bar\phi}(z-v-\eta)_{\mu}^{\mu+\eta{\bar\epsilon}_{k'},j'}
(\sum_{i}
\phi(z-v)_{\mu,i}^{\mu+\eta{\bar\epsilon}_{k'}}
{\bar\phi}(z-v-\eta)_{\mu}^{\mu+\eta{\bar\epsilon}_{k},i})
\phi(z-u)_{\mu,i'}^{\mu+\eta{\bar\epsilon}_{k}}
C_{k}t_{k'}
\eeq
Cancelling by
${\bar\phi}(z-v-\eta)_{\mu}^{\mu+\eta{\bar\epsilon}_{k'},j'}$,
multiplying both sides by
${\bar\phi}(z-v-\eta)_{\lambda}^{\lambda+\eta{\bar\epsilon}_{l},i'}$
and summing over $i'$ we obtain
$$
\sum_{k}(\sum_{j}
\phi(z-u)_{\mu,j}^{\mu+\eta{\bar\epsilon}_{k'}}
{\bar\phi}(z-v-\eta)_{\lambda}^{\lambda+\eta{\bar\epsilon}_{k},j})
(\sum_{i'}
\phi(z-v)_{\lambda,i'}^{\lambda+\eta{\bar\epsilon}_{k}}
{\bar\phi}(z-v-\eta)_{\lambda}^{\lambda+\eta{\bar\epsilon}_{l},i'})
C_{k'}t_{k}=
$$
\beq
=\sum_{k}(\sum_{i}
\phi(z-v)_{\mu,i}^{\mu+\eta{\bar\epsilon}_{k'}}
{\bar\phi}(z-v-\eta)_{\mu}^{\mu+\eta{\bar\epsilon}_{k},i})
(\sum_{i'}
\phi(z-u)_{\mu,i'}^{\mu+\eta{\bar\epsilon}_{k}}
{\bar\phi}(z-v-\eta)_{\lambda}^{\lambda+\eta{\bar\epsilon}_{l},i'})
C_{k}{\tilde t}_{k'}
\eeq
Using formula (\ref{f-f}) and the relation
\beq
\label{vu}
v= u+\sum_{k}(\lambda_{k}-
\mu_{k})
\eeq
we get
$$
\sum_{k}
\Phi_{z-v-\eta}(\lambda_{lk}+\frac{\eta}{n})
\Phi_{z-v-\eta}(\lambda_{k}-\mu_{k'}+\frac{\eta}{n})
t_{k}C_{k'}
\frac{\prod_{m}
\theta(\lambda_{m}-\mu_{k'}+\frac{\eta}{n})
\theta(\lambda_{mk}+\frac{\eta}{n})
}{\prod_{m\neq k}
\theta(\lambda_{mk})
\prod_{s\neq l}
\theta(\lambda_{sl})}=
$$
\beq
\label{phi2}
=\sum_{k}
\Phi_{z-v-\eta}(\lambda_{l}-\mu_{k}+\frac{\eta}{n})
\Phi_{z-v-\eta}(\mu_{kk'}+\frac{\eta}{n})
{\tilde t}_{k'}C_{k}
\frac{\prod_{m}
\theta(\lambda_{m}-\mu_{k}+\frac{\eta}{n})
\theta(\mu_{mk'}+\frac{\eta}{n})
}{\prod_{m\neq k}
\theta(\mu_{mk})
\prod_{s\neq l}
\theta(\lambda_{sl})}
\eeq
At this point we need the following identity (in fact, it was proved in
[NRK]):

  {\bf Lemma.}
$$
\sum_{j}
\Phi_{z}(y_{ij}+\xi)
\Phi_{z}(y_{j}-x_{k}+\xi)
\prod_{m\neq j}
\frac{\theta(y_{jm}-\xi)}
{\theta(y_{jm})}
\prod_{s} \frac{
\theta(x_{s}-y_{j}-\xi)}
{\theta(x_{s}-y_{j})}=
$$
\beq
\label{idphi}
=\sum_{j}
\Phi_{z}(y_{i}-x_{j}+\xi)
\Phi_{z}(x_{jk}+\xi)
\prod_{m\neq j}
\frac{\theta(x_{jm}+\xi)}
{\theta(x_{jm})}
\prod_{s} \frac{
\theta(x_{j}-y_{s}-\xi)}
{\theta(x_{j}-y_{s})}
\eeq
$Proof.$
Using the functional relation for $\Phi_{z}(x)$ (\ref{a}) we find
that it is enought to prove the following identities:
\beq
\label{one}
\sum_{j}
(\zeta(y_{ij}+\xi)+\zeta(y_{j}-x_{k}+\xi))
\prod_{m\neq j}
\frac{\theta(y_{jm}-\xi)}
{\theta(y_{jm})}
\prod_{s}
\frac{\theta(x_{s}-y_{j}-\xi)}
{\theta(x_{s}-y_{j})}=
\eeq
$$
=\sum_{j}
(\zeta(y_{i}-x_{j}+\xi)+\zeta(x_{jk}+\xi))
\prod_{m\neq j}
\frac{\theta(x_{jm}+\xi)}
{\theta(x_{jm})}
\prod_{s}
\frac{\theta(x_{j}-y_{s}-\xi)}
{\theta(x_{j}-y_{s})}
$$
and
\beq
\label{two}
\sum_{j}
\prod_{m\neq j}
\frac{\theta(y_{jm}-\xi)}
{\theta(y_{jm})}
\prod_{s}
\frac{\theta(x_{s}-y_{j}-\xi)}
{\theta(x_{s}-y_{j})}=
\eeq
$$
=\sum_{j}
\prod_{m\neq j}
\frac{\theta(x_{jm}+\xi)}
{\theta(x_{jm})}
\prod_{s}
\frac{\theta(x_{j}-y_{s}-\xi)}
{\theta(x_{j}-y_{s})}.
$$

  From (\ref{b}) we get
\beq
\label{three}
\prod_{l}
\frac{\theta(z-x_{l}+\xi)\theta(z-y_{l}-\xi)}
{\theta(z-x_{l})\theta(z-y_{l})}=
\sum_{j}
(\zeta(z-x_{j})-\zeta(x_{k}-\xi-x_{j}))
\prod_{l}
\frac{\theta(x_{j}-x_{l}+\xi)\theta(x_{j}-y_{l}-\xi)}
{\theta(x_{j}-y_{l})}\times
\eeq
$$
\times
\frac{1}{
\prod_{m\neq j}
\theta(z-x_{l})\theta(x_{j}-y_{l})}+
\sum_{j}
(\zeta(z-y_{j})-\zeta(x_{k}-\xi-y_{j}))\times
$$
$$
\times\prod_{l}
\frac{\theta(y_{j}-x_{l}+\xi)\theta(y_{j}-y_{l}-\xi)}
{\theta(y_{j}-x_{l})}\
\frac{1}{
\prod_{m\neq j}
\theta(y_{j}-y_{m})}.
$$

  Put in this equation $z= y_{i}+\xi$ we get
\beq
\label{four}
\sum_{j}
(\zeta(y_{i}-x_{j}+\xi)-\zeta(x_{k}-\xi-x_{j}))
\prod_{l}
\frac{\theta(x_{j}-x_{l}+\xi)\theta(x_{j}-y_{l}-\xi)}
{\theta(x_{j}-y_{l})}
\frac{1}{
\prod_{m\neq j}
\theta(z-x_{l})\theta(x_{j}-y_{l})}+
\eeq
$$
+\sum_{j}
(\zeta(y_{i}-y_{j}+\xi)-\zeta(x_{k}-\xi-y_{j}))
\prod_{l}
\frac{\theta(y_{j}-x_{l}+\xi)\theta(y_{j}-y_{l}-\xi)}
{\theta(y_{j}-x_{l})}\
\frac{1}{
\prod_{m\neq j}
\theta(y_{j}-y_{m})}= 0.
$$
and cancelling both sides by $\theta(\xi)$ we get (\ref{one}). To prove
(\ref{two}) we use the relation (\ref{c}):
\beq
\label{five}
\sum_{j}
\frac{
\prod_{m}
\theta(y_{j}-x_{m}+\xi)
\theta(y_{j}-y_{m}-\xi)}
{\prod_{m}
\theta(y_{j}-x_{m})
\prod_{m\neq j}
\theta(y_{j}
-y_{m})}+
\eeq
$$
+\sum_{j}
\frac{
\prod_{m}
\theta(x_{j}-x_{m}+\xi)
\theta(x_{j}-y_{m}-\xi)}
{\prod_{m}
\theta(x_{j}-y_{m})
\prod_{m\neq j}
\theta(x_{j}-x_{m})}= 0
$$
and (\ref{two}) follows.
\square

  From the lemma we find that equation (\ref{phi2}) is satisfied then
\beq
\label{nt}
t_{k}= e^{c}\prod_{s}
\frac{\theta(\lambda_{k}-\mu_{s}+\frac{\eta}{n})}
{\theta(\lambda_{k}-\mu_{s})}
\eeq
\begin{equation}
\label{tildet}
{\tilde t}_{k}= e^{c}\prod_{m\neq k}\frac{\theta(\mu_{mk}-\frac{\eta}{n})}
{\theta(\mu_{mk}+\frac{\eta}{n})}
\prod_{s}\frac{\theta(\lambda_{s}-\mu_{k}+\frac{\eta}{n})}
{\theta(\lambda_{s}-\mu_{k})}
\end{equation}
\begin{equation}
\label{nCK}
C_{k}= \prod_{m}\frac{\theta(\mu_{mk}-\frac{\eta}{n})}
{\theta(\lambda_{m}-\mu_{k})}
\end{equation}
i.e. we rederive formulas (\ref{t}) and (\ref{CK}). Also remark that
from our point of view the commutativity of transformations (\ref{nt}),
(\ref{tildet}) is an easy consequence of the commutativity of
elementary modifications at different points.

  Usually, the Lax operator $L(z)_{i}^{j}$ and $M$-operator $M(z)_{i}^{j}$
was written in another gauge. For comparison, we do a corresponding
transformation of our formulas. It is a step from $f$-sheaf of
degree $-1$ to $f$-sheaf of degree zero.
\beq
L(z)_{k',k}=
\sum_{i,j}
\phi(z-v-\eta)_{\lambda,j}^{\lambda+\eta{\bar\epsilon}_{k'}}
{\bar\phi}(z-v-\eta)_{\lambda}^{\lambda+\eta{\bar\epsilon}_{k},i}
L(z)_{i}^{j}= \Phi_{z-v-\eta}(\lambda_{kk'}+\frac{\eta}{n})
\frac{\prod_{l}\theta(\lambda_{lk'}+\frac{\eta}{n})}
{\prod_{l\neq k}\theta(\lambda_{lk})}
t_{k'}
\eeq
\beq
{\tilde L}(z)_{k',k}=
\sum_{i,j}
\phi(z-v-\eta)_{\mu,j}^{\mu+\eta{\bar\epsilon}_{k'}}
{\bar\phi}(z-v-\eta)_{\mu}^{\mu+\eta{\bar\epsilon}_{k},i}
{\tilde L}(z)_{i}^{j}= \Phi_{z-v-\eta}(\mu_{kk'}+\frac{\eta}{n})
\frac{\prod_{l}\theta(\mu_{lk'}+\frac{\eta}{n})}
{\prod_{l\neq k}\theta(\mu_{lk})}
{\tilde t}_{k'}
\eeq
\beq
M(z)_{k',k}=
\sum_{i,j}
\phi(z-v-\eta)_{\mu,j}^{\mu+\eta{\bar\epsilon}_{k'}}
{\bar\phi}(z-v-\eta)_{\lambda}^{\lambda+\eta{\bar\epsilon}_{k},i}
M(z)_{i}^{j}=
\Phi_{z-v-\eta}(\lambda_{k}-\mu_{k'}+\frac{\eta}{n})\times
\eeq
$$
\times\frac{\prod_{l}\theta(\lambda_{l}-\mu_{k'}+\frac{\eta}{n})}
{\prod_{l\neq k}\theta(\lambda_{lk})}
C_{k'}=
\Phi_{z-v-\eta}(\lambda_{k}-\mu_{k'}+\frac{\eta}{n})
\frac{\prod_{l}\theta(\lambda_{l}-\mu_{k'}+\frac{\eta}{n})}
{\prod_{l\neq k}\theta(\lambda_{lk})}
\frac{\prod_{l}\theta(\mu_{lk'}-\frac{\eta}{n})}
{\prod_{i}\theta(\lambda_{l}-\mu_{k'})}
$$
It is obvious that the Lax equation is satisfied. Also note that
$L(z)_{k',k}$ is the map from ${\cal F}_{0}\otimes{\bf\xi}^{-1}$ to
${\cal F}_{0}$, see (\ref{mapff'}).

  As we know, the vector $s_{\mu}(x)=\prod_{l}\theta(x-\mu_{l})$
lies in the kernel of the map
$$
{\cal F}_{-1}|_{P_{u}}
\stackrel{\phi_{u}}{\longrightarrow}
({\cal F}'_{0}\otimes{\bf\xi}^{-1})|_{P_{u}}.
$$
Then we passing to the $f$-sheaf of degree zero, $s_{\mu}(x)$ is mapped
to the column $(s_{\mu}(\lambda_{1}+\frac{\eta}{n}),...,s_{\mu}
(\lambda_{n}+\frac{\eta}{n}))^{t}$ and we must have
\beq
\sum_{k}
L(z)_{k',k}s_{\mu}(\lambda_{k}+\frac{\eta}{n})=
e^{c}s_{\mu}(\lambda_{k}+\frac{\eta}{n})
\eeq
\beq
\sum_{k}
M(z)_{k',k}s_{\mu}(\lambda_{k}+\frac{\eta}{n})=0
\eeq
We show that this is indeed the case. The calculations are parallel to [KS].

  Second equation is a consequence of
$$
\sum_{k}
\frac{\theta(u-v+\lambda_{k}-\mu_{k'}-\frac{n-1}{n}\eta)}
{\theta(\lambda_{k}-\mu_{k'}+\frac{\eta}{n})
\prod_{l\neq k}
\theta(\lambda_{lk})}
\prod_{s}\theta(\lambda_{k}-\mu_{s}+\frac{\eta}{n})=
$$
\beq
=\sum_{k}
\theta(u-v+\lambda_{k}-\mu_{k'}-\frac{n-1}{n}\eta)
\frac{\prod_{s\neq k'}
\theta(\lambda_{k}-\mu_{s}+\frac{\eta}{n})}
{\prod_{l\neq k}
\theta(\lambda_{lk})}= 0
\eeq
where the last equality follows from identity (\ref{c}) and
relation (\ref{vu}),
if we set
$y_{i}=\lambda_{i}, x_{k'}= \mu_{k'}+\frac{n-1}{n}++v-u, x_{i\neq k'}=
\mu_{i}-\frac{\eta}{n}$.

  For the proof of first equation, consider
$$
\sum_{k}
\frac{
\theta(u-v+\lambda_{kk'}-\frac{n-1}{n}\eta)}
{\theta(\lambda_{kk'}+\frac{\eta}{n})}
\frac{\prod_{l}
\theta(\lambda_{lk'}+\frac{\eta}{n})}
{\prod_{l\neq k}
\theta(\lambda_{lk})}
\prod_{s}
\theta(\lambda_{k}-\mu_{s}+\frac{\eta}{n})=
$$
\beq
\label{Ls}
=-\sum_{k}
\theta(v-u+\lambda_{k'k}+\frac{n-1}{n}\eta)
\prod_{s}
\theta(\lambda_{k}-\mu_{s}+\frac{\eta}{n})
\prod_{l\neq k}
\frac{
\theta(\lambda_{lk'}+\frac{\eta}{n})}
{\theta(\lambda_{lk})}.
\eeq

  In [KS] the theta function identity
\beq
\label{KS}
\sum_{k=1}^{n}
\theta(z+x_{k'k}-\xi)
\prod_{s}
\theta(x_{k}-y_{s}+\xi)
\prod_{l\neq k}
\frac{
\theta(x_{k'l}-\xi)}
{\theta(x_{kl})}=
\theta(z)
\prod_{s=1}^{n}
\theta(x_{k'}-y_{s}),
\eeq
where
$$
z= n\xi+\sum_{k}(x_{k}-y_{k})
$$
was established. Putting $\xi= \frac{\eta}{n}, z=\eta+v-u, x_{k}=\lambda_{k},
y_{k}=\mu_{k}$ we get
$$
(\ref{Ls})= -\theta(v-u+\eta)
\prod_{s}
\theta(\lambda_{k'}-\mu_{s}).
$$
Therefore
$$
\sum_{k}
L(z)_{k',k}s_{\mu}(\lambda_{k}+\frac{\eta}{n})=
\frac{t_{k'}}
{\theta(u-v-\eta)}
(-\theta(v-u+\eta))
\prod_{s}
\theta(\lambda_{k'}-\mu_{s})=
$$
\beq
=e^{c}
\prod_{s}
\frac{\theta(\lambda_{k'}-\mu_{s}+\frac{\eta}{n})}
{\theta(\lambda_{k'}-\mu_{s})}
\theta_(\lambda_{k'}-\mu_{s})=
e^{c}s_{\mu}(\lambda_{k'}+\frac{\eta}{n}).
\eeq

  Let us return to formulas (\ref{nt}) and (\ref{tildet}). If we define
the simplectic structure by putting the variables $\lambda_{k},\log t_{k},
k=1,..n$ to be canonically conjugated, then the transformation from
$(\lambda_{k},t_{k})$ to $(\mu_{k},{\tilde t}_{k})$ is canonical. The
generating function of this transformation is
\beq
F(\lambda,\mu)=
\sum_{k,k'}
(S(\lambda_{k}-\mu_{k'}+\frac{\eta}{n})-S(\lambda_{k}-\mu_{k'}))+
\frac{1}{2}
\sum_{k\neq k'}
(S(\mu_{kk'}+\frac{\eta}{n})-S(\mu_{kk'}-\frac{\eta}{n}))+
\eeq
$$
+c(u+\sum_{k}(\lambda_{k}-\mu_{k})),\qquad
S(\lambda)=\int^{\lambda}\log\theta(x)dx
$$
i.e.
\beq
\log t_{k}= \partial F/\partial\lambda_{k},
\log {\tilde t}_{k}= -\partial F/\partial\mu_{k}.
\eeq


  At the end, following
[NRK], we can interpret the equations (\ref{nt}),(\ref{tildet})
as evolution in discrete time. Indeed, put $t_{k}(a)= t_{k},
\lambda_{k}(a)= \lambda_{k}, t_{k}(a+1)= {\tilde t}_{k},
\lambda_{k}(a+1)= \mu_{k}, k=1,...n$ and $a$ is a time index. The
equations (\ref{nt}),(\ref{tildet}) are then
\beq
\label{a}
t_{k}(a)=
e^{c(a)}
\prod_{s}
\frac{\theta (\lambda_{k}(a)-\lambda_{s}(a+1)+\frac{\eta}{n})}
{\theta(\lambda_{k}(a)-\lambda_{s}(a+1))}
\eeq
\beq
\label{a+1}
t_{k}(a+1)=
e^{c(a)}
\prod_{m\neq k}
\frac{\theta(\lambda_{mk}(a+1)-\frac{\eta}{n})}
{\theta(\lambda_{mk}(a+1)+\frac{\eta}{n})}
\prod_{s}
\frac{\theta(\lambda_{s}(a)-\lambda_{k}(a+1)+\frac{\eta}{n})}
{\theta(\lambda_{s}(a)-\lambda_{k}(a+1))}
\eeq

  From (\ref{a+1}) we obtain
\beq
t_{k}(a)=
e^{c(a-1)}
\prod_{m\neq k}
\frac{\theta(\lambda_{mk}(a)-\frac{\eta}{n})}
{\theta(\lambda_{mk}(a)+\frac{\eta}{n})}
\prod_{s}
\frac{\theta(\lambda_{s}(a-1)-\lambda_{k}(a)+\frac{\eta}{n})}
{\theta(\lambda_{s}(a-1)-\lambda_{k}(a))}
\eeq
and combining it with (\ref{a}) we get
\beq
e^{c(a)-c(a-1)}
\prod_{m\neq k}
\frac{\theta(\lambda_{mk}(a)+\frac{\eta}{n})}
{\theta(\lambda_{mk}(a)-\frac{\eta}{n})}=
\prod_{s}
\frac{\theta(\lambda_{k}(a)-\lambda_{s}(a+1))}
{\theta(\lambda_{k}(a)-\lambda_{s}(a+1)+\frac{\eta}{n})}
\frac{\theta(\lambda_{k}(a)-\lambda_{s}(a-1)-\frac{\eta}{n})}
{\theta(\lambda_{k}(a)-\lambda_{s}(a-1))}
\eeq

  This is the discrete  RS equation.

{\centerline {\bf Acknowledgements.}}

  I am indebted to S.Kharchev, A.Marshakov and A.Zabrodin for discussion.

\centerline{\bf References.}

[Bax]R.J.Baxter "Eight-vertex model in lattice statistics and
one-dimensional anizotropic Heisenberg chain" I,II,III Ann.Phys.76(1973)
1-71

[BKMS]V.V.Bazhanov, R.M.Kashaev, V.V.Mangazeev and Yu.G.Stroganov "
$({\bf Z}^{\times}_{N})^{n-1}$ generalizations of the Chiral Potts Model"
CMP 138(1991)393-408

[Dr]V.G.Drinfeld "Varieties of modules of $F$-sheaves" Funct.Analiz i ego
priloz. 21(1987)N2, 23-41

[ER]B.Enriquez, V.Rubtsov "Hitchin systems, higher Gaudin operators
and $r$-matrices" alg-geom/9503010

[Fay]J.Fay "Theta-functions on Riemann surfaces" LNM 352(1973)

[JMO]M.Jimbo,T.Miwa and M.Okado "Solvable lattice models whose
states are dominant integral weights of $A^{(1)}_{n-1}$" Lett.in Math. Phys.
14(1987)123-131

[H1]K.Hasegawa "On the crossing symmetry of the elliptic solution of the
Yang-Baxter equation and a new $L$-operator for Belavin $R$-matrix"
J.Phys.A26(1993)3211

[H2]K.Hasegawa "Ruijsenaars commuting difference operators as commuting
transfer-matrices" CMP187(1997)

[KS]V.B.Kuznetsov and E.K.Sklyanin "On Backlund transformation for
many-body systems" J.Phys.A31(1998)2241-2251

[M]S.Mukai "A duality between $D(E)$ and $D(E^{\Lambda})$ " Nagoya Math. J.
81(1981)153-175

[Mu]D.Mumford "Tata lectures on theta I,II" Birkhauser(1983)

[N]N.Nekrasov "Holomorphic bundles and many body systems" hep-th/9503157

[NRK]F.W.Nijhoff, O.Ragnisco and V.B.Kuznetsov "Integrable
time-discretisation of the Ruijsenaars-Schneider model" CMP176(1996)681-700

[RT]M.P.Richey and C.A.Tracy "${\bf Z}_{n}$ Baxter model: symmetries and
the Belavin parametrization" J.Stat.Phys. 42(1986)311-348

[S]E.K.Sklyanin "Some algebraic structures connected with the
Yang-Baxter equation" Funct. Analiz i ego priloz.16(1982)N4,27-34;
17(1983)N4,34-48.

\end{document}